\documentclass[a4paper, english]{article}

\usepackage{graphicx,indentfirst,amsfonts, amsthm, amsmath, amsthm, amssymb, amscd, amssymb, mathrsfs, makeidx}
\usepackage[all]{xy}

\newtheorem{theorem}{Theorem}[section]
\newtheorem{definition}[theorem]{Definition}

\newtheorem{lemma}[theorem]{Lemma}

\newtheorem{corollary}[theorem]{Corollary}
\newtheorem{proposition}[theorem]{Proposition}

\newenvironment{remark}{\vspace{0.75em} \noindent{\bf Remark.  }}{\vspace{0.75em}}

\newcommand{\C}{\mathbb{C}}

\newcommand{\Z}{\mathbb{Z}}

\newcommand{\Osheaf}{\mathcal{O}}

\newcommand{\Spec}{\mbox{\rm Spec }}

\newcommand{\Pro}{\mathbf{P}}
\newcommand{\pwpr}{\Pro(W) \! \times \Pro^{r}}
\newcommand{\pr}{\Pro^r}

\newcommand{\Hilb}{\mbox{\it{Hilb}}}
\newcommand{\dblq}{/\!/}

\newcommand{\Mgn}{\overline{\mathcal{M}}_{g,n}}

\newcommand{\Hom}{\mbox{\rm Hom}}
\newcommand{\Pic}{\mbox{\rm Pic }}

\newcommand{\acanon}{\omega^{\otimes a}}

\newcommand{\Maps}{\overline{\mathcal{M}}_{g,n}(\mathbf{P}^r,d)}

\newcommand{\pts}{x_1,\ldots ,x_n}
\newcommand{\sections}{\sigma_1,\ldots ,\sigma_n}
\newcommand{\Mapsx}{\overline{\mathcal{M}}_{g,n}(X,\beta)}
\newcommand{\Os}[1]{\mathcal{O}_{#1}}

\newcommand{\PW}[1]{\mathbf{P}(W_{#1})}
\newcommand{\morepts}{x_1,\ldots,x_{n+d}}
\newcommand{\KP}{\mbox{\rm KP}_{\lambda}}

\begin{document}
\title{A GIT Construction of Moduli Spaces of Stable Maps in Positive Characteristic}
\author{Elizabeth Baldwin}
\date{}
\maketitle

\begin{abstract}
In a previous paper, the author and David Swinarski constructed the moduli spaces of stable maps, $\Maps$, via geometric invariant theory (GIT).  That paper required the base field to be the complex numbers, a
restriction which this paper removes: here the coarse moduli spaces of stable maps are constructed via GIT over a more general base.
\end{abstract}

\section*{Introduction}
In \cite{e_and_d}, the author and David Swinarski gave a geometric invariant theory (GIT) construction of $\Maps$, the moduli space of stable maps of degree $d$ from a projective curve with $n$ marked points to projective space $\pr$.  Unfortunately, that construction requires the base to be $\C$ (although the special case of $\Mgn$, the moduli space of stable curves of genus $g$ with $n$ marked points, was constructed over a base field of arbitrary characteristic).  This was because a part of the argument for $\Maps$ depended on the work of Fulton and Pandharipande in \cite{fp}, which is only stated over $\C$.  The methods of \cite{fp} can be extended beyond $\C$, but in any case it is natural to consider a GIT construction over $\Z$, using the work of Seshadri \cite{sesh_completeness, sesh}.  In this paper we show that the construction of \cite{e_and_d} is valid over $\Z[p_1^{-1}\cdots p_j^{-1}]$, where $p_1,\ldots,p_j$ are all primes less than or equal to $d$, by providing a modified proof, which is independent of \cite{fp}.

The construction may not be presented over $\Spec\Z$ itself because the moduli stack of stable maps of degree $d$ fails to be separated over fields characteristic $p$ where $p\leq d$.  Parts of our argument fail in this more complicated situation.

The set-up and language of this paper closely follow \cite{e_and_d}, and are given in Section \ref{set-up}.  In order to present $\Maps$ as a GIT quotient, we first write down a `parameter space of stable maps with extra structure', and then take a quotient to remove the unwanted information.  We follow Fulton and Pandharipande \cite{fp} in choosing this extra structure; an outline is given now.

If $f:(C,\pts)\rightarrow \pr$ is a stable map then there is a canonically defined line bundle $\mathcal{L}$ on $C$; see line (\ref{L}) below.  For $a\geq 3$ the bundle $\mathcal{L}^a$ is ample.  We fix $a\geq 5$ and subject to further constraints (discussed after Definition \ref{j}).  Now a choice of basis of $H^0(C,\mathcal{L}^a)$ is the extra data we associate to our stable map.

Such a choice of basis provides an embedding $C\hookrightarrow \PW{}$, where $W$ is a fixed vector space of dimension $h^0(C,\mathcal{L}^a)$ (note that this integer depends only on the genus $g$, degree $d$ and number of markings $n$ of the stable map).  It is now logical to take, for our `parameter space with extra structure', a subscheme of the Hilbert scheme of $n$-pointed genus $g$ curves in $\pwpr$;  the subscheme $J$ which we take is that which models curves embedded in $\PW{}$ via $\mathcal{L}^a$, and such that the projection morphism $p_r:(C,\pts)\rightarrow \pr$ is a stable map.  The group $SL(W)$ acts on $J$, such that the orbits are precisely the isomorphism classes of stable maps (Proposition \ref{properties of j}).  It follows that a geometric quotient of $J$ will be the moduli space $\Maps$.  

We wish to form a GIT quotient of a projective scheme, and so we let $\bar{J}$ be the closure of $J$ in the Hilbert scheme.  Now, if a linearisation $L$ of the action of $SL(W)$ on $\bar{J}$ satisfies $\bar{J}^{ss}(L)=\bar{J}^s(L)=J$, then it follows that $\bar{J}\dblq_L SL(W)$ is a geometric quotient of $J$, and so is equal to $\Maps$.

In \cite{e_and_d}, we find a range of linearisations such that this is the case, although some of the proofs in that paper are only valid over $\C$.  The main theorem of this paper is that the construction works more generally; as laid out in Section \ref{main theorem}, the GIT quotient $\bar{J}\dblq_L SL(W)$ is indeed equal to $\Maps$ for such $L$, over any base over which the moduli stack of stable maps is separated.  

To prove this we replace a part of the argument of \cite{e_and_d}.  In Section \ref{main theorem} we give the statements of the required extra results, as well as stating the theorems which we shall use from \cite{e_and_d}.  The results we use from \cite{e_and_d}, which are valid in general characteristic, are that for a range of linearisations of the action of $SL(W)$ on $\bar{J}$, the semistable set is non-empty and contained in $J$; to be precise, the non-emptiness is proved by induction (see Theorem \ref{inductive} for more details).  

The further ingredients we need to prove Theorem \ref{theorem} are that the semistable set $\bar{J}^{ss}$ is an open and closed subset of $J$ (shown in Section \ref{closed section}), and that the scheme $J$ is connected (shown in Section \ref{connected section}).  Along with the knowledge that the semistable set is non-empty, this is sufficient to show that it is, as required, the whole of $J$.

I would like to thank Frances Kirwan for all her help, and Dan Abramovich and Christian Liedtke for useful comments on the first version of this paper, as well as the Engineering and Physical Sciences Research Council and the Mittag-Leffler Institute for funding and supporting this research.
\section{Definitions and Notation}\label{set-up}
Our aim is to construct as a GIT quotient the moduli space $\Maps$ of stable maps of degree $d$ from an $n$-pointed genus $g$ curve.  Fundamental to this construction is the following result, which we quote from Newstead:
\begin{proposition}[\cite{new} Proposition 2.13]
\label{propertiesoffamily}\label{newstead}
Suppose that the family $\mathcal{X}\rightarrow S$ has the local universal property for some moduli problem, and that the algebraic group $G$ acts on $S$, with the property that $\mathcal{X}_s \sim \mathcal{X}_t$ if and only if $G\cdot s = G\cdot t$.  Then:
\begin{enumerate}
\item[(i)]
any coarse moduli space is a categorical quotient of $S$ by $G$;
\item[(ii)]
a categorical quotient of $S$ by $G$ is a coarse moduli space if and only if it is an orbit space.\hfill $\Box$ 
\end{enumerate}
\end{proposition}

Fix a projective scheme $X$.  The discrete invariant $\beta$ may be thought of as a class in $H_2(X;\Z)$; to be precise we look to \cite{behrend-manin} Definition 2.1, and define
\[ 
H_2(X)^+ = \{\alpha\in\Hom_{\Z}(\Pic X, \Z) \ | \ \alpha(L) \geq 0 \mbox{ whenever L is ample}\}.
\]
If $C$ is a connected curve and $f:C\rightarrow X$ is a morphism, then $L\mapsto \deg f^* L$ defines an element of $H_2(X)^+$ which we shall call the homology class of $f$, and denote by $f_*[C]$.  Now we may define our moduli problem:
\begin{definition}
A \emph{stable map} in $\Mapsx$ is a
map $f:(C,\pts) \rightarrow X$ where $(C,\pts)$ is an $n$-pointed prestable curve of genus $g$, the homology class $f_*[C]=\beta$, and the following stability conditions are satisfied: if $C'$ is a nonsingular rational component of $C$ and $C'$ is mapped to a point by $f$, then $C'$ must have at least three {\emph special points} (either marked points or nodes); if $C'$ is a component of arithmetic genus $1$ and $C'$ is mapped to a point by $f$, then $C'$ must contain at least one special point.  
\end{definition}
Note that since we require the domain curves $C$ to be connected, the stability condition on genus one components is automatically satisfied except in $\overline{\mathcal{M}}_{1,0}(X,0)$, which is empty.

If $(C,\pts)$ is a connected reduced projective nodal curve of genus $g$, with $n$ distinct non-singular marked points $\pts$, then a morphism $f:(C,\pts)\rightarrow\pr$ is a stable map if and only if the invertible sheaf 
\begin{equation}
\mathcal{L}:=\omega_C(x_1+\cdots +x_n)\otimes f^*\Osheaf_{\Pro^r}(c),\label{L}
\end{equation}
is ample, where $c$ is any integer greater than or equal to $3$.  In fact greater generality is possible; unless $(g,n,d)=(0,0,1)$, we require only $c\geq 2$.  A discussion on the magnitude of $c$ may be found in \cite{e_and_d} Section 2.4.

Let $a$ be an integer such that $a\geq 5$; then $\mathcal{L}^a$ is very ample (further comments on the magnitude of $a$ are given in the remark below).  Let 
\[
e:=\deg \mathcal{L}^a=a(2g-2+n+cd)
\]
and let $W=W_{g,n,d}$ be a vector space of dimension
\[
N+1:=e-g+1=h^0(C,\mathcal{L}^a).
\]
We will not use the subscripts $g,n,d$ on $W$ except where they are needed for clarity.  A choice of basis for $H^0(C,\mathcal{L}^a)$ determines an embedding $C\rightarrow\PW{}$; taking the graph of $f$ we obtain $C$ as a subscheme in $\pwpr$, such that the map $f$ is retrieved by projecting to $\pr$.  Moreover, $C\subset\pwpr$ satisfies 
\begin{equation}\label{sheaf condit}
(\mathcal{O}_{\Pro(W)}(1) \otimes \mathcal{O}_{\Pro^{r}}(1))|_{C}\cong (\omega_{C}^{a}(a x_1+\cdots +a x_n) \otimes \mathcal{O}_{\Pro^{r}}(c a+1))|_{C}.
\end{equation}

On the other hand, if we have a reduced nodal connected $n$-pointed curve $(C,\pts)\subset\pwpr$, where the marked points are non-singular and distinct, and such that (\ref{sheaf condit}) is satisfied, we conclude that the morphism $p_r:C\rightarrow \pr$ given by projection to $\pr$ is a stable map.  Note that we always use $p_r$ to denote the projection morphism $\pwpr\rightarrow \pr$, or the restriction of this morphism to any subscheme $X\subset\pwpr$.

These ideas motivate our choice of the scheme from which we shall form a quotient.  It is a subscheme of the Hilbert scheme of $n$-pointed curves in $\pwpr$.  To define that, we write $\Hilb(\pwpr)$ for the Hilbert scheme of genus $g$ curves in $\pwpr$ with bidegree $(e,d)$.
\begin{definition}[\cite{fp} Section 2.3]\label{i}\index{$I$}
The scheme 
\[
I=I_{g,n,d}\subset \Hilb(\pwpr)\times\prod_{i=1}^n(\pwpr) 
\]
is the closed incidence subscheme consisting of those $n+1$-tuples $(h,\pts)$ in $\Hilb(\pwpr)\times(\pwpr)^{\times n}$ such that the points $x_1,\ldots,x_n$ lie on $\mathcal{C}_h$.  
\end{definition}

The scheme $I$ is thus the Hilbert scheme of $n$-pointed curves of bidegree $(e,d)$ in $\pwpr$.  The universal family $\mathcal{C}\rightarrow Hilb(\pwpr)$ gives rise to a universal family $(\mathcal{C}^I\rightarrow I, \sections)$ of $n$-pointed curves in $\pwpr$, possessing $n$ sections $\sections:I\rightarrow \mathcal{C}^I$ (which give the marked points).  We wish to consider the subscheme of $I$ corresponding to $a$-canonically embedded stable maps:
\begin{definition}[\cite{fp} Section 2.3]\label{j}\index{$J$}
The scheme $J=J_{g,n,d}\subset I_{g,n,d}$ is the locally closed subscheme consisting of those $(h,\pts)\in I$ such that:
\begin{enumerate}
\item[(i)] $\mathcal{C}_h$ is projective, connected, reduced and nodal, and $\pts$ are non-singular, distinct points on $\mathcal{C}_h$;
\item[(ii)] the projection map $\mathcal{C}_{h} \rightarrow \Pro(W)$ is a non-degenerate embedding;
\item[(iii)] $ (\mathcal{O}_{\Pro(W)}(1) \otimes \mathcal{O}_{\Pro^{r}}(1))|_{\mathcal{C}_{h}} $  and $(\omega_{\mathcal{C}_{h}}^{a}(a x_1+\cdots +a x_n) \otimes \mathcal{O}_{\Pro^{r}}(c a+1))|_{\mathcal{C}_{h}} $ are isomorphic.
\end{enumerate}
We denote by $\bar{J}=\bar{J}_{g,n,d}$ the closure of $J$ in $I$.
\end{definition}
\noindent As with $W$, we only specify the genus, number of markings and degree associated to $I$ and $J$ where it is necessary to do so for clarity.

\begin{remark} Though it is only necessary for $a$ to be at least $3$ to ensure that $\mathcal{L}^a$ is very ample, this is not adequate to provide the GIT construction of $\Maps$: it is shown in \cite{Schubert} that cusps are stable if $a=3$.  The bound $a\geq 10$ is sufficient to ensure that the GIT quotient is isomorphic to $\Maps$ for all values of $g,n,d$, though this is unnecessarily prescriptive for all but the smallest values of $g,n,d$, when $a\geq 5$ will suffice.  Precise details are given in the statement of \cite{e_and_d} Theorem 6.1.
\end{remark}

The Hilbert scheme is defined over $\Spec \Z$, and so it follows that $I$ and $J$ are also.  However, we cannot complete the constructions of this paper over the entirety of $\Z$, as the moduli stack is not separated in this generality (see \cite{behrend-manin} 4.1 and 4.2).  For simplicity of notation we consider the schemes $I$ and $J$ to be defined over $\Z[p_1^{-1}\cdots p_j^{-1}]$, where $p_1,\ldots,p_j$ are all the prime numbers less than or equal to the degree $d$ of the maps in question.  If $k$ is a field then we shall denote $J_{g,n,d}\times_{\Spec \Z}\Spec k$ by $(J_{g,n,d})_k$, or by $J_k$ if in the context it is not necessary to be more specific.

We restrict the family $\mathcal{C}^I\rightarrow I$ over $J$, denoting the restriction $\mathcal{C}^J$; with projection to $\pr$ we obtain the family:
\[
\xymatrix{ \mathcal{C}^J \ar@{^{(}->}@<-0.5ex>[r] \ar[d] & \pwpr\times J \ar[dl]^{p_3} \ar[r]^(0.7){p_r} & \pr \\
		J \ar@<1ex>[u]^{\sigma_i} 	&&}
\]
where $p_r$ is projection to the factor $\pr$.  The family $(\mathcal{C}^J\rightarrow J,\sections,p_r)$ is seen to be a family of stable maps, since as discussed above the map $p_r:(C,\pts)\rightarrow \pr$ is stable precisely when the sheaf $(\omega_{C}^{a}(a x_1+\cdots +a x_n) \otimes \mathcal{O}_{\Pro^{r}}(c a))|_{C}$ is very ample.

There is a natural action of $SL(W)$ on $\PW{}$.  We define $SL(W)$ to act trivially on $\pr$, and so induce an action of the group $SL(W)$ on the Hilbert scheme $\Hilb(\pwpr)\times(\pwpr)^{\times n}$.  The subschemes $I$ and $J$ are easily seen to be invariant under this action.

The significance of the scheme $J$ for moduli of stable maps is now summarised in the following proposition:
\begin{proposition}[\cite{e_and_d} Proposition 3.4]\label{properties of j}
\begin{enumerate}
\item[(i)]
$(\mathcal{C}^J\rightarrow J,\sections,p_r)$ has the
local universal property for the moduli problem of stable maps, $\Maps$.
\item[(ii)]$(h,x_1,\ldots ,x_n)\in J$ and $(h',x'_1 ,\ldots x'_n)\in J$ lie in the same orbit under the action of $SL(W)$ if and only if the stable maps $(\mathcal{C}_{h},x_1,\ldots ,x_n,p_r|_{\mathcal{C}_h})$ and $(\mathcal{C}_{h'},x'_1,\ldots ,x'_np_r|_{\mathcal{C}_{h'}})$ are isomorphic.\hfill$\Box$
\end{enumerate}
\end{proposition}
\noindent It follows, by Proposition \ref{newstead}, that a geometric quotient of $J$ by $SL(W)$ will be the coarse moduli space $\Maps$.

We wish to construct $\Maps$ as a GIT quotient of the projective subscheme $\bar{J}$ of $I$ by $SL(W)$.  Such a quotient is a categorical quotient of its semistable set, and this quotient is in addition a geometric quotient if all points in the semistable set are in fact stable.  Thus we will know that
\[
\bar{J}^{ss}\dblq_L SL(W) \cong \Maps
\]
if we can find a linearisation $L$ such that
\begin{equation}
\bar{J}^{ss}(L)=\bar{J}^{s}(L)=J.\label{needed}
\end{equation}
The principles of GIT were laid out by Mumford \cite{git}, and extended over $\Z$ by Seshadri \cite{sesh_completeness, sesh}.

In \cite{e_and_d} we find a linearisation satisfying (\ref{needed}); however, some of our work there depends on the independent construction of $\Maps$ by Fulton and Pandharipande in \cite{fp}, and is only valid over $\C$.  This is unsatisfactory; in this paper we will remove the dependence on independent constructions and on base $\C$.  

It is somewhat complicated to define the linearisations of the action of $SL(W)$ on $\bar{J}$ that we use, and we shall not need to know their precise form to prove the results of this paper.  Details are given in \cite{e_and_d} Section 4; the fact that Theorems \ref{subset}, \ref{nonempty} and \ref{inductive} hold is all that is needed for this paper.


\section{The Results and structure of this paper}\label{results from thesis}\label{main theorem}

This paper extends the results of \cite{e_and_d} to positive characteristic, by proving (\ref{needed}) over a more general base.  The main theorem of this paper is the following:
\begin{theorem}\label{theorem}
Let $L$ be a linearisation of the action of $SL(W)$ on $\bar{J}$.
\begin{enumerate}
\item[(i)]Let $k$ be a field of characteristic $0$ or $p$ where $p>d$.  Suppose that
\begin{equation}
\emptyset\neq \bar{J}_k^{ss}(L)\subset J_k. \label{sufficient in theorem}
\end{equation}
Then
\begin{equation*}
\bar{J}_k^s(L)=\bar{J}_k^{ss}(L)=J_k 
\end{equation*}
and so
\[
\bar{J}_k\dblq_L SL(W)\cong (\Maps)_k.
\]
In particular, if (\ref{sufficient in theorem}) holds for all such fields, then
\[
\bar{J}\dblq_L SL(W)\cong \Maps
\]
as schemes over $\Z[p_1^{-1}\cdots p_j^{-1}]$, where $p_1,\ldots,p_j$ are the prime numbers less than or equal to $d$.
\item[(ii)]Suppose that $L$ satisfies the conditions of \cite{e_and_d} Theorem 6.1.  The same conclusions hold:
\[
\bar{J}^s(L)=\bar{J}^{ss}(L)=J 
\]
and so
\[
\bar{J}\dblq_L SL(W)\cong (\Maps).
\]
as schemes over $\Z[p_1^{-1}\cdots p_j^{-1}]$, where $p_1,\ldots,p_j$ are the prime numbers less than or equal to $d$.
\end{enumerate}
\end{theorem}
\begin{remark}We emphasise that, for any choice of $g$, $n$ and $d$ giving rise to stable maps, there do exist linearisations satisfying the conditions of \cite{e_and_d} Theorem 6.1.
\end{remark}

As a corollary we will be able to construct the moduli space of stable maps to a general projective scheme $X$.
\begin{corollary}[\cite{e_thesis} Corollary 3.2.8, cf.\ \cite{fp} Lemma 8]
Let $X$ be a projective scheme, with a fixed embedding to projective space $X\stackrel{\iota}{\hookrightarrow}\Pro^r$, and let $\beta\in H_2(X)^+$.  Let $g$ and $n$ be non-negative integers.  If $\beta=0$ then suppose in addition that $2g-2+n\geq 1$.  Let $\iota_*(\beta)=d\in H_2(\Pro^r)^+$.

Then there exists a closed subscheme $J_{X,\beta}$\index{$J_{X,\beta}$} of $J_{g,n,d}$, and there exist linearisations $L$ of the action of $SL(W)$ (namely any $L$ such that the equality $\bar{J}_{g,n,d}^{ss}(L)=\bar{J}_{g,n,d}^s(L)=J_{g,n,d}$ is satisfied and in particular any linearisation satisfying the conditions of \cite{e_and_d} Theorem 6.1), with
\[
\bar{J}_{X,\beta} \dblq_{L|_{\bar{J}_{X,\beta}}} SL(W) \cong \overline{\mathcal{M}}_{g,n}(X,\beta),
\]
where $\bar{J}_{X,\beta}$ is the closure of $J_{X,\beta}$ in $\bar{J}$.  Again this is an isomorphism of schemes over $\Z[p_1^{-1}\cdots p_j^{-1}]$, where $p_1,\ldots,p_j$ are the prime numbers less than or equal to $d$.\hfill$\Box$
\end{corollary}
The majority of the proofs given in \cite{e_and_d} are valid over an arbitrary base.  We state here the ones which we shall use to prove Theorem \ref{theorem}.
\begin{theorem}[\cite{e_and_d} 5.21, see also \cite{e_thesis} 5.6.1]\label{subset}
Let $k$ be any field of any characteristic.
If $L$ is a linearisation of the action of $SL(W)$ on $\bar{J}$ satisfying the hypotheses given in \cite{e_and_d} Theorem 6.1, then
\[
\bar{J}_k^{ss}(L)\subset J_k.
\]
\hfill$\Box$
\end{theorem}
\begin{theorem}[\cite{e_and_d} 6.5, cf.\ \cite{gieseker} 1.0.0]\label{nonempty}
Let $k$ be any field of any characteristic.  Suppose that a smooth curve $C\subset\pwpr$, defined over $k$, is represented in $(J_{g,0,d})_k$.  Suppose $L$ is a linearisation of the action of $SL(W)$ on $(\bar{J}_{g,0,d})_k$, which satisfies the hypotheses given in \cite{e_and_d} Theorem 6.1.  Then $C$ is semistable in $(\bar{J}_{g,0,d})_k$ with respect to $L$.\hfill$\Box$
\end{theorem}
\begin{theorem}[\cite{e_and_d} 6.9, see also \cite{e_thesis} 6.3.1]\label{inductive}
Let $k$ be any field, and suppose that $n>0$.  Let $L_{g+1,n-1,d}$ be a linearisation of the action of $SL(W_{g+1,n-1,d})$ on $(\bar{J}_{g+1,n-1,d})_k$, satisfying the hypotheses of \cite{e_and_d} Theorem 6.1, as they apply to stable $(n-1)$-pointed maps of genus $g+1$ and degree $d$.  Suppose moreover that 
\[
(\bar{J}_{g+1,n-1,d})_k^{ss}(L_{g+1,n-1,d})=(\bar{J}_{g+1,n-1,d})_k^{s}(L_{g+1,n-1,d})=(J_{g+1,n-1,d})_k
\]
Then, for any linearisation $L_{g,n,d}$ satisfying the hypotheses of \cite{e_and_d} Theorem 6.1 applied to stable $n$-pointed maps of genus $g$ and degree $d$, we have:
\[
(\bar{J}_{g,n,d})_k^{ss}(L_{g,n,d})=(\bar{J}_{g,n,d})_k^{s}(L_{g,n,d})=(J_{g,n,d})_k.
\]
\hfill$\Box$
\end{theorem}

Note that Theorem \ref{nonempty}, whose proof closely follows that of \cite{gieseker} Theorem 1.0.0, is only valid for $n=0$ (though in that case it does in fact hold for a wider range of linearisations than is stated here).  The method of proof fails to show semistability for smooth maps with $n>0$ marked points, with respect to any linearisations that satisfy the hypotheses of \cite{e_and_d} Theorem 6.1.

We use the  `inductive step' provided by Theorem \ref{inductive}, as Theorem \ref{nonempty} is not sufficient on its own to give semistability of any stable maps with marked points, with respect to the key linearisations.  It should be clarified that \cite{e_and_d} provides suitable linearisations $L_{g,n,d}$ for the action of $SL(W_{g,n,d})$ on the space $\bar{J}_{g,n,d}$ for each combination of $g,n,d$ giving rise to stable maps.  For fixed values of $g,n,d$, Theorem \ref{inductive} relates  semistability of points in the scheme $\bar{J}_{g,n,d}$ with respect to suitable linearisations $L_{g,n,d}$, with semistability of points in the different space $\bar{J}_{g+1,n-1,d}$, with respect to suitable linearisations $L_{g+1,n-1,d}$ for this space.

It follows from Theorem \ref{inductive} that, if we can prove (\ref{needed}) for the base case $n=0$ for all genera $g$, over a given field $k$, then (\ref{needed}) holds over $k$ for all $n$.  It remains, then, to show that this base case holds, in as great a generality as possible.  

The new technical results that we will prove in this paper are the following:
\begin{proposition}\label{closed}
Suppose $k$ has characteristic $0$ or $p$ where $p>d$.  If $L$ is a linearisation of the action of $SL(W)$ on $\bar{J}_k$ such that $\bar{J}_k^{ss}(L)\subset J_k$, then $\bar{J}_k^{ss}(L)$ is closed as a subscheme of $J_k$.
\end{proposition}
\begin{proposition}\label{j_connected}
Suppose that $k$ is an algebraically closed field.  Then $J_{k}$ is connected.
\end{proposition}
The proof of Proposition \ref{closed} is given in Section \ref{closed section}, and that of Proposition \ref{j_connected} is given in Section \ref{connected section}.
\begin{proof}[Proof of Theorem \ref{theorem}, assuming Propositions \ref{closed} and \ref{j_connected}]
(i) Assume that 
\[
\emptyset\neq \bar{J}_k^{ss}(L)\subset J_k.
\]
Then the stabiliser subgroup of every semistable point is zero-dimensional, and so $\bar{J}_k^{ss}(L)=\bar{J}_k^{s}(L)$ by \cite{git} page 10.  

By definition the semistable set is an open subscheme of $\bar{J}_k$, and hence it is open in $J_k$.  By Proposition \ref{closed}, the semistable set $\bar{J}_k^{ss}(L)$ is in addition a closed subscheme of $J_k$, since the characteristic of $k$ is $0$ or $p$ where $p>d$.  We thus know that $\bar{J}_k^{ss}(L)$ is a non-empty open and closed subset of $J_k$.

On the other hand, Proposition \ref{j_connected} says that $J_{\bar{k}}$ is connected, where $\bar{k}$ is the algebraic closure of $k$.  It follows that the three schemes in (\ref{needed}) are the same over $\bar{k}$, i.e.\ 
\[
\bar{J}_{\bar{k}}^{ss}(L)=\bar{J}_{\bar{k}}^{s}(L)=J_{\bar{k}}.
\]
However, the stable and semistable sets of a scheme are defined by their geometric points (\cite{git} Definition 1.7).  The geometric points $X_k(\bar{k})$ of a scheme $X_k$ over $k$ are in bijection with $X_{\bar{k}}(\bar{k})$ (cf.\ \cite{liu} 2.18(a)).  Thus equality (\ref{needed}) follows.

Now suppose (\ref{sufficient in theorem}) holds over any field $k$ of characteristic $0$ or $p$ where $p>d$.  We show that the two open subsets $\bar{J}^{ss}(L)\subset\bar{J}$ and $J\subset\bar{J}$ are equal as schemes over $\Z[p_1^{-1}\cdots p_j^{-1}]$ (where $\bar{J}^{ss}(L)$ is defined over $\Z[p_1^{-1}\cdots p_j^{-1}]$ as in \cite{sesh_completeness, sesh}).  As $\bar{J}^{ss}(L)$ and $J$ are both \emph{open} subschemes of $\bar{J}$, it is sufficient to show that they are set-theoretically equal, as then it automatically follows that the scheme structure is the same.

Let $x$ be any point in $\bar{J}^{ss}(L)$; we may consider $x$ as a point in $\bar{J}_{k(x)}^{ss}(L)$, where $k(x)$ is the residue field at $x$.  We know that $\bar{J}_{k(x)}^{ss}(L)=J_{k(x)}$, and so $x$ lies also in $J_{k(x)}$ and hence in $J$; thus $\bar{J}^{ss}(L)\subset J$.  Similarly we show $J\subset \bar{J}^{ss}(L)$, and hence $\bar{J}^{ss}(L)=J$ over $\Z[p_1^{-1}\cdots p_j^{-1}]$.  The same argument may be used to show that these schemes are equal to the semistable set, $\bar{J}^{s}(L)$, and so we obtain the required equality 
\[
\bar{J}^{ss}(L)=\bar{J}^s(L)=J
\]
of schemes over $\Z[p_1^{-1}\cdots p_j^{-1}]$, recalling again that $p_1,\ldots,p_j$ are the prime numbers less than or equal to $d$.

(ii) By Theorem \ref{inductive} and the argument given at the end of part (i), it is sufficient to show that for every field $k$ of suitable characteristic,
\begin{equation}
(\bar{J}_{g,0,d})_k^{ss}(L_{g,0,d})=(\bar{J}_{g,0,d})_k^{s}(L_{g,0,d})=(J_{g,0,d})_k  \label{base case needed}
\end{equation}
holds, for any linearisation $L_{g,0,d}$ of the group action on $\bar{J}_{g,0,d}$, such that $L_{g,0,d}$ satisfies those conditions of \cite{e_and_d} Theorem 6.1 which pertain to genus $g$ degree $d$ stable maps with no marked points.  

However, by Theorems \ref{subset} and \ref{nonempty} we know that
\[
\emptyset\neq(\bar{J}_{g,0,d})_k^{ss}(L_{g,0,d})\subset (J_{g,0,d})_k, 
\]
so (\ref{base case needed}) follows by part (i), as required.
\end{proof}

It remains, then, to prove Propositions \ref{closed} and \ref{j_connected}, so we know that the semistable set is closed in $J_k$, and that $J_k$ is connected, over suitable fields $k$.  These results are given in the following sections.

\begin{remark}The question of what happens when the characteristic of $k$ is less than or equal to $d$ remains.  Theorems \ref{subset} and \ref{nonempty} still hold, so when $n=0$ we obtain a projective quotient which is a geometric quotient of some nonempty open subscheme of $J$.  This then coarsely represents some substack of the moduli stack of stable maps, which in this case is non-separated.
\end{remark}

\section{The semistable set is closed in $J$}\label{closed section}
Here we use semistable replacement, together with the fact that the moduli stack of stable maps is separated, to prove Proposition \ref{closed}: that if $L$ is a linearisation of the action of $SL(W)$ on $\bar{J}_k$ such that the semistable set $\bar{J}_k^{ss}(L)\subset J_k$, then $\bar{J}_k^{ss}(L)$ is a closed subscheme of $J_k$, provided that $k$ has characteristic $0$ or $p$ where $p>d$.  The moduli stack is only proved to be separated over such base fields; the bound on the characteristic ensures that the stable maps themselves are separable morphisms.  For more details see \cite{behrend-manin} Proposition 4.1 and Lemma 4.2.  This is why the bound on the characteristic of $k$ is needed in Theorem \ref{theorem}.

\begin{proof}[Proof of Proposition \ref{closed}]
We shall show that the inclusion $\bar{J}_k^{ss}(L)\hookrightarrow J_k$ is a proper morphism, using the valuative criterion of properness.  As we work with a fixed suitable linearisation $L$, we shall abbreviate the notation to $\bar{J}^{ss}_k$.
Let $R$ be a discrete valuation ring with quotient field $K$.  Suppose we have:
\[
\xymatrix{
\Spec K\ar[r] \ar@{^{(}->}[d]	 	& \bar{J}_k^{ss} \ar@{^{(}->}@<-1ex>[d]\\
\Spec R \ar[r]^(0.6)f				& J_k }.
\]
The proposition is proved if we can complete the diagram with a morphism from $\Spec R$ to $\bar{J}_k^{ss}$.  It is sufficient to show that $f(0_R)\in \bar{J}_k^{ss}$, where $0_R$ is the unique closed point in $R$.

By the Cohen Structure Theorem (\cite{hartshorne} Theorem 5.5A), the completion of $R$ is a formal power series ring $A$ over a field $L$.
Then we have:
\[
\xymatrix{
\Spec L \ar [r] \ar@{^{(}->}[d] 	&  \Spec K \ar[r] \ar@{^{(}->}[d]	 	& \bar{J}_k^{ss} \ar@{^{(}->}[d]\\
\Spec A \ar[r]				& \Spec R \ar[r]^(0.6)f				& J_k }.
\]
Now, by \cite{sesh_completeness} Theorem 4.1, there exists an element $g\in SL(W)(K)$ and a finite cover $\Spec A \rightarrow \Spec A$ such that we obtain a morphism $\beta':\Spec A\rightarrow \bar{J}^{ss}$ and the following diagram commutes:
\[
\xymatrix{
\Spec L \ar [r] \ar@{^{(}->}[d] & \Spec K \ar[r] \ar@{^{(}->}[d]	& \bar{J}_k^{ss} \ar@{^{(}->}[d] \ar[r]^{\cdot g} 	
											&\bar{J}_k^{ss} \ar@{^{(}->}[r] 		& J_k\\
\Spec A \ar[r]	\ar@{.>}[urrr]^{\beta'} \ar@/_{1pc}/[rr]_{\alpha}	& \Spec R \ar[r]			& J_k							&& }.
\]
By composition, we now have two morphisms $\alpha, \beta:\Spec A \rightarrow J_k$, where $\alpha$ composes our original morphisms along the lower level and $\beta$ composes $\beta'$ with the inclusion $\bar{J}_k^{ss}\hookrightarrow J_k$.  The morphisms $\alpha$ and $\beta$ are in general distinct.  However, the images $\alpha(\Spec L)$ and $\beta(\Spec L)$ in $J_k$ of the generic point $\Spec L$ of $\Spec A$ are in the same orbit under the action of $SL(W)$; they are related by our fixed element $g\in SL(W)(K)$.  It follows that the stable maps corresponding to these two points of $J_k$ are isomorphic, since points in the same orbit under the action of $SL(W)$ give rise to isomorphic stable maps (Proposition \ref{properties of j}).

We shall show that the two images $\alpha(0_A)$ and $\beta(0_A)$ in $J_k$ of the closed point $0_A$ in $\Spec A$ are also in the same orbit as one another.  By definition, $\beta(0_A)$ is in $\bar{J}_k^{ss}$, which is $SL(W)$-invariant, so it will follow that $f(0_R)=\alpha(0_A)\in\bar{J}_k^{ss}$, as required.

We may pull back the universal family $(\mathcal{C}^{J_k}\rightarrow J_k, \sections,p_r)$ along $\alpha$ and $\beta$ to obtain two families over $\Spec A$ of stable $n$-pointed maps from curves embedded in $\Pro(W)$.  The two families are isomorphic on the generic point $\Spec L\subset \Spec A$.  We now use the fact that the moduli functor for stable maps is separated (for $d$ bounded by characteristic), shown in \cite{behrend-manin} Proposition 4.1\@.  It follows that the two families are isomorphic families of stable maps over $\Spec A$.

In particular, the stable maps given by the fibre of either family over $0_A$ are isomorphic, and hence $\alpha(0_A)$ and $\beta(0_A)$ lie in the same orbit in $J_k$ under the action of $SL(W)$.  It follows, since $\beta(0_A)\in \bar{J}_k^{ss}$, and since $\bar{J}_k^{ss}$ is a $SL(W)$-invariant subscheme of $J_k$, that $f(0_R)=\alpha(0_A)\in\bar{J}_k^{ss}$, as required.
\end{proof}

\section{The scheme $J_k$ is connected}\label{connected section}

It remains to prove Proposition \ref{j_connected}: that the scheme $J_k$ is connected, provided that $k$ is an algebraically closed field.  Our argument is very similar to that of Kim and Pandharipande in \cite{kp}.
They construct a family of a special type of stable maps, over a base scheme which is connected, and show that the image of this base in $\Mapsx$ must meet every connected component of $\Mapsx$.  However, to apply this argument to $J$ we must construct a family of stable maps whose domains are consistently embedded in $\PW{}$ via $\mathcal{L}^a$.  It is embedding the entire family into $\pwpr$ which requires the extra work.  

On the other hand, as we prove this result only for maps to projective space, we are able to work over a more general base than Kim and Pandharipande; the method requires use of the Kleiman-Bertini theorem, which is only proved over the complex numbers in general, but over any algebraically closed field for the special case of projective space (see \cite{kleiman-transversality} Corollary 11).

We work over a fixed algebraically closed base field $k$ throughout this section, so we shall omit the subscripts $k$ throughout; other subscripts are constantly in play and the notation would become too cumbersome.

We fix a 1-PS $\lambda:\mathbb{G}_m\rightarrow SL_{r+1}$, which acts on $\pr$ with distinct weights.  We shall show that, generically, stable maps of a special type emerge as the limit points of the action of $\lambda$ on $\pr$.  First, if $p\in \pr$ is a point fixed by $\lambda$, we may define
\begin{eqnarray*}
A_p	&:=& \{x\in\pr | \lim_{t\rightarrow 0}\lambda(t)\cdot x = p \},\\
A'_p	&:=& \{x\in\pr | \lim_{t\rightarrow \infty}\lambda(t)\cdot x = p \}.
\end{eqnarray*}
Let $0\in\pr$ be the unique fixed point such that $A_0$ is open in $\pr$; then $A'_0=\{0\}$.  Let $1\in\pr$ be the unique fixed point such that $A_1$ is codimension $1$ in $\pr$, so that $A'_1$ is one-dimensional.  Define 
\[
U:=A_0\cup A_1.  
\]
This is an open subscheme whose complement in $\pr$ is codimension $2$.  Set 
\[
P:=\overline{A'_1}=A'_1\cup A'_0.  
\]
Then $P$ is isomorphic to $\Pro^1$.  We also note that $P\subset U$, and that $P\cap \overline{A_1}=\{1\}$.  The easiest way to verify all these statements is to choose homogeneous coordinates for $\pr$, with respect to which the action of $\lambda$ is diagonal.  
In addition, $P$ is a generator for the homology group $H_2(\pr)^+$ (as defined in Section \ref{set-up}). The line $P$ is associated to the positive generator of homology given by the morphism $L\mapsto \deg(L|_P)$.  

For convenience, we introduce a name for the special maps distinguished by Kim and Pandharipande:
\begin{definition}[\cite{kp} Section 2]\label{def_KP}
A stable map $f:(C,\pts)\rightarrow \pr$ of genus $g$ is \emph{of type $\KP$} if it satisfies the following conditions:
\begin{enumerate}
\item[(i)] $C$ is of the form $\tilde{C}\cup\bigcup_{i=1}^d \Pro^1_i$, where $\tilde{C}$ is a prestable curve of genus $g$ and each $\Pro^1_i$ is a projective line meeting $\tilde{C}$ in a single node, at $y_i$; 
\item[(ii)] all the markings $\pts$ lie on $\tilde{C}$;
\item[(iii)] $f(\tilde{C})=0\in\pr$ and $f|_{\Pro^1_i}:\Pro^1_i\cong P$ for $i=1,\ldots,d$.
\end{enumerate}
\end{definition}
\noindent Note that this definition depends on the specific choice of $\lambda$ (which defines $0$ and $P$).  Since $f$ is a stable map which collapses the subcurve $\tilde{C}$ to a point, it follows that $(\tilde{C},\pts,y_1,\ldots,y_d)$ is a stable  $(n+d)$-pointed curve.
\begin{proposition}[cf.\ \cite{kp} Proposition 2]\label{kp_in_all}
Suppose $k$ is algebraically closed.  Every connected component of $J$ contains a stable map of type $\KP$\@.
\end{proposition}
\begin{proof}
When one works in the moduli space of stable maps, one may show (\cite{kp} Proposition 2) that the stabilisation of the limit point of a generic stable map under the action of $\lambda$ on $\pr$ is a map of type $\KP$.  It only takes a little more work to extend these ideas to $J$.

Recall our universal family of $n$-pointed stable maps over $J$:
\[
\xymatrix{ \mathcal{C}^J \ar@<-0.4ex>@{^{(}->}[r]	\ar[d]	& \pwpr	\times J \ar[dl]^{p_3} \ar[r]^(.7){F}	& \pr \\
	 J \ar@<1ex>[u]^{\sigma_i}			& &}
\]
where $\sections$ are the sections corresponding to the marked points, and $F$ is projection onto the second factor.
There is an action of $SL_{r+1}$ on $J$, induced by the action on $\pr$.  
Any two points in the same orbit under the action of $SL_{r+1}$ will lie in the same connected component of $J$, as $SL_{r+1}$ itself is connected.

Let $J'$ be a connected component of $J$ and let $(h,\pts)\in J'$.  By the Kleiman-Bertini Theorem (\cite{kleiman-transversality} Corollary 11) we know that a general $SL_{r+1}$-translate 
\[
f=g\cdot F_h:(C,\pts)\rightarrow\pr 
\]
of the stable map $F_h$ satisfies the following properties (cf.\ \cite{kp} Section 2):
\begin{enumerate}
\item[(i)] the image $f(C)$ lies in $U$, i.e.\ does not meet the codimension $2$ complement of $U$;
\item[(ii)] $f(C)$ intersects the codimension $1$ subscheme $A_1$ transversally at non-singular points of $C$;
\item[(iii)] the markings of $C$ have image in the open subscheme $A_0$.
\end{enumerate}
Since $f$ and $F_h$ lie in the same orbit of $J$ under the action of $SL_{r+1}$, it follows that $f$ is also represented in the connected component $J'$ of $J$.

Forget the embedding of $C$ in $\Pro(W)$, and refer to the stable map from this abstract curve as $[f]$.  Return to our fixed 1-PS $\lambda:\mathbb{G}_m\rightarrow SL_{r+1}$.  We use this action to obtain a limit stable map which is of type $\KP$\@.  Define a family of stable $n$-pointed maps over $\mathbb{G}_m$:
\[
\xymatrix{\mathbb{G}_m \times C \ar[r]^(.55){h} \ar@<0.5ex>[d] 	&\pr\\
		\mathbb{G}_m \ar@<0.5ex>[u]^{\sigma_i} &}
\]
where each section $\sigma_i$ is defined by the marked point $x_i\in C$, and the morphism $h$ is given by the action of $\lambda$ on $f(C)\subset\pr$.  We embed $\mathbb{G}_m$ in $\mathbb{A}^1_k$.  As the moduli stack of stable maps is proper (\cite{behrend-manin} Theorem 3.14) our family extends over $\mathbb{A}^1_k$.  The fibre over $t=1\in\mathbb{G}_m$ is precisely the stable map $[f]$.

Details of this construction may be found in \cite{kp} Proposition 2, where it is proved the extension of our family to a family 
\[
(S\rightarrow \mathbb{A}_k^1,\sections,h) 
\]
over $\mathbb{A}_k^1\supset \mathbb{G}_m$ of stable maps to $\Pro^1$ has the property that the fibre of $S$ over the limit point $t=0\in\mathbb{A}_k^1$ is of type $\KP$\@.  We refer to this limit stable map as $[h_0]$.  We shall say a few word on why intuitively this is so.

Consider the limit $\lambda(t)\cdot f$ as $t$ tends to $0$.  One may picture the points of $f(C)$ which lie in $A_0$ being contracted to $0$, while the points which meet $A_1$ are contracted to $1$, leaving rational tails between $0$ and $1$, namely on $P$.  Such a map is a good candidate to be a stable map of type $\KP$.  In fact this limit map is not in general a stable map; we may end up with rational components with fewer than three special points in the limit curve, which are collapsed by the limit map.  However, we may `stabilise' our family by contracting such rational components to a point.

As stated in Proposition \ref{properties of j}, the scheme $J$ possesses the local universal property for the moduli problem of stable maps.  Thus there exists an open neighbourhood $U$ of $0\in\mathbb{A}_k$ and a morphism $\phi:U\rightarrow J$ which induces a family of stable maps isomorphic to $S|_U\rightarrow U$.  

The point $1\in\mathbb{A}_k$ is not necessarily also in $U$, but we can in any case find a finite connected open covering $\mathcal{U}=\{U_i\}$ of $[0,1]\subset\mathbb{A}_k^1$ and a collection of morphisms $\{\phi_i:U_i\rightarrow J\}$, each inducing an family of stable maps isomorphic to the corresponding $S|_{U_i}$.  Moreover, the images $\phi_i(U_i)$ of these open subsets all lie in the same connected component of $J$\@.  For if $y\in U_i\cap U_j\neq\emptyset$ then $\phi_i(y)$ and $\phi_j(y)$ are both points in $J$ representing the stable map $h_y:S_y\rightarrow\pr$ and so by Proposition \ref{properties of j} they lie in the same orbit of $J$ under the action of $SL(W)$.  As each $U_i$ is itself connected, we may conclude that we have found representatives for $[h_0]$ and $[h_1]=[f]$ which are in the same connected component of $J$.

Recall that $f$ lies in $J'$, and so the orbit of $f$ in $J$ under the action of $SL(W)$ is contained in $J'$.  Our new representative for $[f]$ lies in this orbit and so in $J'$; it follows that our representative for $[h_0]$ does also.  This is a map of type $\KP$, and so the proof is completed.
\end{proof}


We now show that all stable maps of type $\KP$ are represented in the same connected component of $J$.  For this, we must work with spaces corresponding to maps of varying genera, number of marked points and degree, thus we must not abbreviate the notation $J_{g,n,d}$ and $W_{g,n,d}$.

\begin{proposition}\label{family_exists}
There exists a family:
\[
\xymatrix{ \mathcal{X}	\ar@<-0.4ex>@{^{(}->}[r] \ar[d]	& \Pro(W_{g,n,d})\times\pr\times B \ar[r]^(0.75){p_2} \ar[dl]^{p_3}	& \pr\\
	B \ar@<1ex>[u]^{\sigma_i} &&}
\]
of $n$-pointed stable maps of type $\KP$, modelling all such stable maps, and such that $B$ is connected.
\end{proposition}
\begin{remark}The fibres of the family we shall construct in this proof shall in fact satisfy the defining conditions for $J$; see Proposition \ref{family in j}.
\end{remark}
\begin{proof}
Recall that if $f:(C,\pts)\rightarrow\pr$ is a stable map of type $\KP$, with genus $g$ component $\tilde{C}$ meeting the rational components at nodes $y_1,\ldots,y_d$, then $(\tilde{C},\pts,y_1,\ldots y_d)$ is a stable $(n+d)$-pointed curve of genus $g$.  Conversely, such a curve defines an unique stable map of type $\KP$; we attach projective lines at the last $d$ marked points and define the map as described in \ref{def_KP}.  However, to prove this proposition we must embed the domain curve in $\Pro(W_{g,n,d})$, in a way which is consistent over a family.  

We shall create our new family by modifying the universal family of stable $a$-canonically embedded curves in $\Pro(W_{g,n+d,0})$:
\[
\xymatrix{\mathcal{C}^{J_{g,n+d,0}} \ar[d] \ar@<-0.4ex>@{^{(}->}[r] 	& \PW{g,n+d,0}\times J_{g,n+d,0} \ar[dl]^{p_3} \\
	J_{g,n+d,0} \ar@<1ex>[u]^{\sigma_i}.			&}  
\]
First, we need to be able to embed any curve represented in $J_{g,n+d,0}$ into $\Pro(W_{g,n,d})\times\pr$, so we fix an embedding $\Pro(W_{g,n+d,0})\subset\Pro(W_{g,n,d})$ as follows.

Let $V_{0,1,1}$ be a vector space of dimension $a(c-1)-1=\dim W_{0,1,1}-1$.  By comparing dimensions, we see that we may choose and fix a decomposition
\begin{equation}
W_{g,n,d}\cong W_{g,n+d,0}\oplus \bigoplus_{i=1}^d V_{0,1,1}^i \label{main_decomp}.
\end{equation}
where $V_{0,1,1}^1,\ldots, V_{0,1,1}^d$ are distinguished copies of $V_{0,1,1}$.  Thus we may project from $W_{g,n,d}$ to any one of the additive factors, and so view $\Pro(W_{g,n+d,0})$ and $\Pro(V_{0,1,1}^i)$ for $i=1,\ldots,d$ as fixed linear subspaces of $\Pro(W_{g,n,d})$.  

Finally, define a map $\mathcal{C}^{J_{g,n+d,0}}\rightarrow \pr$ collapsing the whole family to the specified point $0$; now we have an embedding
\begin{equation}
\mathcal{C}^{J_{g,n+d,0}}\hookrightarrow \PW{g,n,d}\times\pr\times J_{g,n+d,0}.\label{embedded C}
\end{equation}
We shall first describe the glueing process for a single curve, and use this to induce the whole new family. Suppose that 
\[
(\tilde{C},\morepts)\subset\Pro(W_{g,n+d,0})\times\{0\}\subset\Pro(W_{g,n,d})\times\pr 
\]
is a stable curve represented in $J_{g,n+d,0}$.  Note in particular that projection to $\pr$ collapses $\tilde{C}$ to the point $0\in\pr$.  Now let 
\[
(D,y)\subset\Pro(W_{0,1,1})\times\pr 
\]
be a smooth genus zero curve, represented in $J_{0,1,1}$, such that $p_r(D)=P$ and $p_r(y)=0$, where $P$ and $0$ are our fixed line and point in $\pr$. We shall attach d copies of $(D,y)$ to $\tilde{C}$ by identifying each $y$ sequentially with $x_{n+1},\ldots,x_{n+d}$.

Having identified $W_{0,1,1}$ with $H^0(\Pro(W_{0,1,1}),\Osheaf_{\Pro(W_{0,1,1})}(1))$, we note that this is of dimension $a(c-1)$ and so we may fix an isomorphism between the codimension $1$ subspace of sections vanishing at $y$, and the vector space $V_{0,1,1}$.  Choose some section $u\in W_{0,1,1}$ which is non-zero at $y$; we obtain a decomposition, which we shall fix:
\begin{equation}
W_{0,1,1}\cong\langle u\rangle \oplus V_{0,1,1}. \label{split_up_011}
\end{equation}

On the other hand, we let $V_{g,n+d,0}$ be the subspace of $W_{g,n+d,0}$ consisting of sections vanishing at all the marked points $x_{n+1},\ldots,x_{n+d}$, and choose $d$ independent sections 
\[
u_1,\ldots,u_d\in H^0(\Pro(W_{g,n,d}),\Os{\Pro(W_{g,n,d})}(1)),
\]
such that $u_i(x_{n+j})=\delta_{ij}$, for $i,j=1,\ldots,d$.  We may identify any one of these with $u\in W_{0,1,1}$, and so may set $W_{0,1,1}^i:=\langle u_i \rangle \oplus V_{0,1,1}^i$ to be a distinguished copy of $W_{0,1,1}$, for $i=1,\ldots,d$.  Now the decomposition (\ref{main_decomp}) induces another decomposition:
\begin{equation}
W_{g,n,d}=V_{g,n+d,0}\oplus\bigoplus_{i=0}^d W_{0,1,1}^i. \label{example_decomp}
\end{equation}
Thus we obtain projections $W_{g,n,d}\rightarrow W_{0,1,1}^i$ for $i=1,\ldots, d$, and so identify $d$ corresponding linear subspaces $\Pro(W_{0,1,1}^i)\subset\Pro(W_{g,n,d})$; note that these are pairwise disjoint.  We had a genus zero curve $(D,y)\subset\Pro(W_{0,1,1})\times P$, and so a copy 
\[
(D_0^i,y_0^i)\subset \Pro(W_{0,1,1}^i)\times P \subset \PW{g,n,d}\times\pr
\]
for $i=1,\ldots d$.  The subscript $0$ reflects the fact that we wish to fix these as curves with specific embeddings in $\PW{g,n,d}\times\pr$. 
 
 All sections in $V_{0,1,1}^i$ vanish at $y_0^i$, as do all sections in $W_{0,1,1}^j$ for $i\neq j$, and all sections in $V_{g,n+d,0}$.  In addition, we specified that $p_r(y_0^i)=0\in\pr$, so it follows that the image of $y_0^i$ in $\Pro(W_{g,n,d})\times\pr$ is at precisely the point $x_{n+i}\in \tilde{C}$.  Thus $D_0^i$ and $\tilde{C}$ meet here, and this singular point is a node as the curves lie in two linear subspaces meeting transversally.  On the other hand, none of the rational components meet one another, as they live in disjoint linear subspaces.  We obtain a new $n$-pointed curve, 
\[
(C,\pts):=(\tilde{C}\cup D^1\cup\cdots\cup D^d,\pts)\subset\Pro(W_{g,n,d})\times\pr,
\]
where the unions are taken of subschemes of $\PW{g,n,d}\times\pr$.  It is clear by construction that $(C,\pts)\stackrel{p_r}{\rightarrow}\pr$ is a stable map of type $\KP$.  

When we extend this to the universal family $\mathcal{C}^{J_{g,n+d,0}}\rightarrow J_{g,n+d,0}$, the point at which we attach each copy of $(D,y)$ must vary; we vary the way in which we embed them into $\Pro(W_{g,n,d})\times\pr$, as follows.

The group $G:=SL(W_{g,n+d,0})$ acts on $\Pro(W_{g,n+d,0})$; if we define it to act trivially on the other factors in our direct sum (\ref{main_decomp}), and trivially on $\pr$, then we have an action of $G$ on $\Pro(W_{g,n,d})\times\pr$ which leaves the linear subspace $\Pro(W_{g,n+d,0})\times\pr$ invariant.  There is a non-trivial action on the sections $u_1,\ldots,u_d$, and so the action is non-trivial on the embedded curves $D^1_0,\ldots,D^d_0$.

Now we may induce $d$ families of $1$-pointed curves over $G$; they are all trivial but come with distinct non-trivial embeddings in $\Pro(W_{g,n,d})\times\pr$.  The product $G\times D_0^i$ is embedded in $\PW{g,n,d}\times\pr\times G$ via the action map of $G$, so we have, for $i=1,\ldots, d$, a family:
\[
\xymatrix{ 
\mathcal{D}^i:=G\times D_0^i\, \ar[d] \ar@<-0.4ex>@{^{(}->}[r]
		& \Pro(W_{g,n,d})\times\pr\times G \ar[dl]^{p_3}\ar[r]^(0.55){p_{1,2}}	& \Pro(W_{g,n,d})\times\pr\\
G \ar@<1ex>[u]^{\gamma^i} \ar@<-1ex>[urr]_{\pi^i}.		&					&}
\]
The section $\gamma^i$ is given by $g\mapsto (g,y_0^i)$ for $g\in G$, and so the marked point in the fibre over $g$ is at $(g\cdot y_0^i,g)\in\PW{g,n,d}\times\pr\times G$.  The morphism $\pi^i$ is the composition of this with projection to $\PW{g,n,d}\times\pr$, and so is defined by
\[
\pi^i:g\mapsto g\cdot y_0^i
\]
and its significance is explained in the following.

We intend to take a fibre product of these families with the universal family $\mathcal{C}^{J_{g,n+d,0}}\rightarrow J_{g,n+d,0}$, attaching the curves in the fibres at the images of the sections $\gamma^i$ and $\sigma_{n+i}$.  The morphism $\pi^i$ maps a point $g$ in $G$, the base of the family, to the location of $\gamma^i(g)\in\PW{g,n,d}\times\pr$ in the corresponding fibre.  We shall define similar morphisms from $J_{g,n+d,0}$ to $\PW{g,n,d}\times\pr$, and then use these to define our fibre product.

As the $\pi^i$ are significant, then, we analyse them further.  Recall that $y_0^i\in D_0^i$ lies in the linear subspace $\PW{g,n+d,0}\times\{0\}$ of $\PW{g,n,d}\times\pr$.  This subspace is invariant under the action of $G$, and so $\pi^i$ maps $G=SL(W_{g,n+d,0})$ into $\PW{g,n+d,0}\times\{0\}\cong \PW{g,n+d,0}$.  Let $H\subset G$ be the subgroup scheme stabilising $y_0^i$.  Then $H$ acts on $G$ on the right.  By \cite{git} Theorem 1.1 and Amplification 1.3, a universal geometric quotient $(Y,\phi)$ of $G$ by $H$ exists; in particular the morphism $\phi:G\rightarrow Y$ is universally submersive and surjective.  However it is clear that  $\pi^i:G\rightarrow \Pro(W_{g,n+d,0})$ is in fact this quotient morphism, and so it follows that, for $i=1,\ldots,d$, the morphism $\pi^i$ is surjective, universally submersive and has fibres isomorphic to $H$.  Moreover $H$ may readily be seen to be connected, so $\pi^i$ has connected fibres.

We wish to attach $d$ rational curves to each curve in the family $\mathcal{C}^{J_{g,n+d,0}}\rightarrow J_{g,n+d,0}$ in such a way that the rational curves are disjoint from one another.  Let us then note where each fibre of $\mathcal{D}^i\rightarrow G$ lies in $\PW{g,n,d}\times\pr$.  The curve we started with, $D_0^i$, lies in the linear subspace spanned by the fixed subspace $\Pro(V_{0,1,1}^i)$, and by $y_0^i\in\PW{g,n+d,0}\times\pr$; hence for each $g\in G$, the fibre $\mathcal{D}^i_g$ lies in the linear subspace spanned by $\Pro(V_{0,1,1}^i)$ and by $g\cdot y_0^i$.  In particular, the only point at which $\mathcal{D}^i_g$ meets the linear subspace $\PW{g,n+d,0}\times\pr$ is at $g\cdot y_0^i$. 

The other family we use in our fibre product construction is the universal family over $J_{g,n+d,0}$, embedded in $\PW{g,n,d}\times\pr$ as at line (\ref{embedded C}), namely:
\[
\xymatrix{ \mathcal{C}^{J_{g,n+d,0}} \ar@{^{(}->}[r]\ar[d] 
		& \PW{g,n,d}\times\pr\times J_{g,n+d,0} \ar[dl]^{p_3}\ar[r]^(0.6){p_{1,2}} & \PW{g,n,d}\times\pr \\
	J_{g,n+d,0}. \ar@<1ex>[u]^{\sigma_i} \ar@<-1ex>[urr]_{s_i}  &&}
\]
We have disjoint sections $\sigma_i:J_{g,n+d,0}\rightarrow\mathcal{C}^{J_{g,n+d,0}}$ for $i=1,\ldots,n+d$.  For the final $d$ marked points, however, we are more interested in the maps
\[
s_i:=p_{1,2}\circ\sigma_{n+i}:J_{g,n+d,0}\rightarrow\PW{g,n+d,0}\times\pr
\]
for $i=1,\ldots d$.

We shall attach the rational components at the marked points one at a time, and inductively show that we obtain a family as desired. The base case of this inductive process is the family
\[
\mathcal{X}^0:=\mathcal{C}^{J_{g,n+d,0}}\rightarrow B^0:=J_{g,n+d,0},
\]
with sections $\sigma_i^0:B^0\rightarrow\mathcal{X}^0$ for $i=1,\ldots,n+d$ and morphisms $s_i^0:B^0\rightarrow\PW{g,n,d}\times\pr$ for $i=1,\ldots,d$.

We make inductive definitions, for $j=1,\ldots,d$.
\begin{eqnarray*}
B^{j}	&:=& B^{j-1} {}_{s_{j}^{j-1}}\!\!\times_{\pi^{j}} G \\
	&&	\mbox{ with projections $p^j_{1}:B^j\rightarrow B^{j-1}$ and $p^j_{2}:B^j\rightarrow G$}. \\
\mathcal{X}^j &:=& p^{j\,*}_1\mathcal{X}^{j-1} \cup p^{j\,*}_2 \mathcal{D}^j, \\
		&&	\mbox{ where the union is of two subschemes of $\PW{g,n,d}\times\pr\times B^j$}.
\end{eqnarray*}
We pull the sections $\sigma_i^{j-1}$ back over $p^j_{1}$ to obtain sections $\sigma_i^j:B^j\rightarrow\mathcal{X}^j$, for $i=1,\ldots,n+d$.  We let $s_i^j:B^j\rightarrow \PW{g,n,d}\times\pr$ be given by $p_{1,2}\circ\sigma_{n+i}^j$, for $i=1,\ldots,d$.
The picture is:
\begin{equation}
\xymatrix{\mathcal{X}^j:=p_1^*\mathcal{X}^{j-1}\cup p_2^*\mathcal{D}_{j}	\ar[rr]\ar[dr]\ar[dd] 	&& \mathcal{D}^{j} \ar[dr] 	& \\
		&	B^j \ar[dd]^{p^j_{1}}\ar[rr]^{p^j_{2}}\ar@<1ex>[ul]^{\sigma_i^{j}}	 && G\ar[dd]^{\pi^j} \\
	\mathcal{X}^{j-1}\ar[dr] 			&&				& \\					
		&	B^{j-1}\ar[rr]^{s_j^{j-1}}\ar@<1ex>[ul]^{\sigma_i^{j-1}} && \PW{g,n,d}\times\pr. }\label{cube}
\end{equation}
To complete the proof of Proposition \ref{family_exists} it remains to show:
\begin{lemma}\label{ind_step}
For $j=0,\ldots,d$, the family $(\mathcal{X}^j\rightarrow B^j,\sigma_1^j,\ldots,\sigma_{n+d}^j)$ defined above is a family of curves in $\PW{g,n,d}\times\pr$ with $n+d$ distinct marked points, 
satisfying:
\begin{enumerate}
\item for each $b\in B^j$ the fibre $\mathcal{X}^j_b$ is a nodal curve, consisting of a genus $g$ subcurve, $\tilde{\mathcal{C}}_b$, which has $n+d$ marked points, $\sigma_i^j(b)$ for $i=1,\ldots,n+d$, and lies in the linear subspace $\PW{g,n+d,0}\times\{0\}$, and of $j$ rational components $\mathcal{D}^i_b$, for $i=1,\ldots,j$, each meeting $\tilde{\mathcal{C}}_b$ at the marked point $\sigma^j_{n+i}(b)$, and each lying in the subspace of $\PW{g,n,d}\times\pr$ spanned by the corresponding $\Pro(V_{0,1,1}^i)\times P$ and by $\sigma_{n+i}^{j}(b)$;
\item the base $B^j$ is connected.
\end{enumerate}
\end{lemma}
\emph{Proof of Lemma \ref{ind_step}.} We check that this is true for the base case $j=0$.  Part 1 we know by construction.  For 2 we note that the moduli space of curves $\overline{\mathcal{M}}_{g,n+d}$ is well-known to be connected, and that there is a surjective, universally submersive morphism $J_{g,n+d,0}\rightarrow\overline{\mathcal{M}}_{g,n+d}$ whose fibres are isomorphic to $SL(W_{g,n+d,0})$, witnessed in \cite{e_and_d} Theorem 6.3 (this construction of that paper is valid over arbitrary base field); hence $B^0=J_{g,n+d,0}$ is connected. 

Make the inductive hypothesis that the lemma holds for the family $\mathcal{X}^{j-1}\rightarrow B^{j-1}$.  We check the two parts in turn.

1. For any $b\in B^{j}$, the fibres of $p^{j\,*}_1\mathcal{X}^{j-1}$ and $p^{j\,*}_2 \mathcal{D}^j$ over $b$ meet transversally in $\PW{g,n,d}\times\pr$, at $s_{j}^{j-1}(b)=\pi^{j}(b)$.  The fibres over $b$ are otherwise disjoint; for $1\leq i\leq j-1$ we inductively know where each rational component $\mathcal{D}^i_b$ lies, and $p^{j\,*}_2\mathcal{D}^j_b$ lies, by construction in the linear subspace spanned by $\Pro(V_{0,1,1}^j)$ and by $\sigma_{n+j}^j(b)$; this meets the genus $g$ component only at $\sigma^j_{n+j}(b)=s_j^j(b)$, and none of the other rational components, since the sections $\sigma_i^j$ are disjoint.  It follows that the family $\mathcal{X}^j\rightarrow B^j$ is nodal, each component lying in the linear subspace described.

2. Connectedness is a property of the underlying topological spaces of our schemes.  We shall use the fact that, if $f:X\rightarrow Y$ is a surjective submersive morphism of topological spaces such that $Y$ is connected and such that all the fibres of $f$ are connected, then $X$ is connected.

We have assumed $B^{j-1}$ to be connected, and have constructed $B^j$ via the Cartesian square
\[
\xymatrix{ B^j \ar[r]^{p^j_{2}} \ar[d]^{p^j_{1}}	& G \ar[d]^{\pi^j}\\
	B^{j-1}\ar[r]^{s_j^{j-1}}			&\PW{g,n,d}\times\pr. }
\]
The fibres of $p^j_{1}$ are the fibres of $\pi^j$, which we have already seen to be connected.  Moreover the morphism $\pi^j$ is universally submersive and surjective, and so it follows that $p^j_{1}$ is surjective and submersive.  

This completes the inductive step, and proves the lemma; when we forget the final $d$ sections, the resulting family completes the proof of Proposition \ref{family_exists}.
\end{proof}
Moreover, we may now prove:
\begin{proposition}\label{family in j}
Every fibre $(\mathcal{X}^b,\sigma_i(b),\ldots\sigma_n(b))\subset\pwpr$ of the family constructed in the proof of Proposition \ref{family_exists} satisfies conditions (i)-(iii) in the definition of $J_{g,n,d}$.
\end{proposition}
\begin{proof}
We are concerned with the family given by the $d$th step of the induction; as we no longer need the preceding families in the induction, we may drop the superscript $d$ and so denote it:
\[
\xymatrix{ \mathcal{X}\, \ar[d] \ar@<-0.4ex>@{^{(}->}[r] & \pwpr\times B \ar[dl]^{p_3} \ar[r]^(0.75){p_2} & \pr \\
		B. \ar@<1ex>[u]^{\sigma_i} &&}
\]
By construction, every fibre gives a stable map of type $\KP$, whose domain is non-degenerately embedded in $\PW{g,n,d}$.  Thus conditions (i) and (ii) in the definition \ref{j} of $J$ are verified.  It remains to show that each fibre satisfies \ref{j}(iii), the sheaf condition. 

In general, if $C$ is a nodal curve and $C'\subset C$ is a complete subcurve, meeting the rest of $C$ in only one node at $Q$, then
\[
\omega_C|_{C'}= \omega_{C'}(Q).
\]
If we look at the curve in our family over $b\in B$, we may write
\[
\mathcal{X}_d=\tilde{C}\cup D^1\cup\cdots\cup D^d
\]
where $\tilde{C}$ has genus $g$ and each of the $D^i$ is rational, for $i=1,\ldots,d$.  Let $x_1,\ldots,x_n$ be the marked points (which all lie on $\tilde{C}$) and, for $i=1,\ldots,d$ let $y_i$ be the node where $\tilde{C}_d$ meets $D^i$.  By the construction of $\mathcal{X}\rightarrow B$, the $(n+d)$-pointed curve $(\tilde{C}_b,\pts,y_1,\ldots,y_d)\subset\PW{g,n+d,0}\times\{0\}\subset\PW{g,n,d}\times\pr$ is modelled in $J_{g,n+d,0}$, and so:
\begin{eqnarray}
\lefteqn{ \left(\Osheaf_{\Pro(W_{g,n,d})}(1)\otimes\Osheaf_{\Pro^r}(1)\right)|_{\tilde{C}}
	= \Osheaf_{\Pro(W_{g,n+d,0})}(1)|_{\tilde{C}} }\notag\\
	&\hspace{1.5in}&\cong \acanon_{\tilde{C}}(ax_1+\cdots+ax_{n}+ay_1+\cdots+ay_d) \label{tilde c cong}\\
		&& = \left(\acanon_{\mathcal{X}_d}(ax_1+\cdots+ax_n)\otimes \Osheaf_{\Pro^r}(ca+1)|_{\mathcal{X}_d}\right)|_{\tilde{C}}. \notag
\end{eqnarray}
On the other hand, for every $i=1,\ldots,d$, each pointed rational curve $(D^i,y_i)$ is a linear translation of the fixed embedded curve $(D_0^i,y_0^i)\subset\PW{0,1,1}\times P \subset \PW{g,n,d}\times\pr$, and that embedded curve $(D_0^i,y_0^i)\subset\PW{0,1,1}\times\pr$ is represented in $J_{0,1,1}$; hence the following is satisfied, for $i=1,\ldots,d$:
\begin{eqnarray}
\lefteqn{ \left(\Osheaf_{\Pro(W_{g,n,d})}(1)\otimes\Osheaf_{\Pro^r}(1)\right)|_{D^i}
	= \left(\Os{\Pro(W_{1,1,0})}(1)\times\Os{\pr}(1)\right)|_{D^i} }\notag\\
	&\hspace{1.5in}&\cong \acanon_{D^i}(ay_i)\otimes\Os{\pr}(ca+1)|_{D^i} \label{Di cong}\\
		&&= \left(\acanon_{\mathcal{X}_d}(ax_1+\cdots+ax_n)\otimes\Osheaf_{\Pro^r}(ca+1)|_{\mathcal{X}_d}\right)|_{D^i} .\notag
\end{eqnarray}
To extend these isomorphisms of line bundles over the whole of $\mathcal{X}_d$, we simply need to insist that the isomorphisms over $\tilde{C}$ and $D^i$ are consistent at $y_i$, for $i=1,\ldots,d$.
When we restrict to the fibre over $y_i$, the two isomorphisms (\ref{tilde c cong}) and (\ref{Di cong}) are scalar multiples of one another, so we obtain consistency at $y_i$ for $i=1,\ldots,d$ by multiplying each isomorphism (\ref{Di cong}) by a suitable non-zero scalar, once (\ref{tilde c cong}) is given.

It follows that the fibre $(\mathcal{X}_d,\pts)$ satisfies condition (iii) in the definition of $J$, which completes the proof of Proposition \ref{family in j}.
\end{proof}
Now we put these pieces together, to obtain:
\begin{proposition}\label{KP all together}
All points in $J_{g,n,d}$ representing a stable map of type $\KP$ are in the same connected component of $J_{g,n,d}$.
\end{proposition}
\begin{proof}
$\mathcal{X}\rightarrow B$ is a family of $n$-pointed curves in $\Pro(W_{g,n,d})\times\Pro^r$, so by the universal property there exists a morphism $\Phi:B\rightarrow I_{g,n,d}$.  As each fibre satisfies conditions (i)-(iii) in the definition of $J_{g,n,d}\subset I_{g,n,d}$, it follows that the image of $\Phi$ lies in $J_{g,n,d}$.  Since $B$ is connected it follows that $\Phi(B)$ is connected.

Moreover, $\mathcal{X}\rightarrow B$ represents all stable maps of type $\KP$.  If any point in $J_{g,n,d}$ represents a stable map of type $\KP$, then this map is also represented in $\Phi(B)$.  These two representatives in $J_{g,n,d}$ must therefore lie in the same orbit of the action of $SL(W_{g,n,d})$ on $J_{g,n,d}$, and hence lie in the same connected component.
\end{proof}

\begin{proof}[Proof of Proposition \ref{j_connected}]
The connectedness of $J$ now follows from Propositions \ref{kp_in_all} and \ref{KP all together}.  
\end{proof}
This completes the proofs of the technical results needed for Theorem \ref{theorem}.

\end{document}